# *Convex Relaxation for Optimal Fixture Layout Design*

Zhen Zhong[a], Shancong Mou[a*], Jeffrey H. Hunt[b], and Jianjun Shi[a]

[a]*H. Milton Stewart School of Industrial and Systems Engineering, Georgia Institute of Technology, Atlanta, the United States;*

[b]*The Boeing Company, El Segundo, the United States*

**Abstract**

This paper proposes a general fixture layout design framework that directly integrates the system equation with the convex relaxation method. Note that the optimal fixture design problem is a large-scale combinatorial optimization problem, we relax it to a convex semidefinite programming (SDP) problem by adopting sparse learning and SDP relaxation techniques. It can be solved efficiently by existing convex optimization algorithms and thus generates a near-optimal fixture layout. A real case study in the half-to-half fuselage assembly process indicates the superiority of our proposed algorithm compared to the current industry practice and state-of-art methods.

## 1. Introduction

Composite materials are widely used in the aircraft industry due to their superior properties including high strength-to-weight ratio, corrosion resistance, and high durability (Jones 1998). Aircraft fuselages consist of thin sheets of composite materials of a large area which makes them compliant parts (Megson 2016). Due to the compliant property, the deformation and dimensional variation induced by gravity load will significantly influence the quality of the assembly process. When the deformation is larger than the engineering specification, a shape adjustment method has to be adopted to adjust the shape of the fuselage, which will not only introduce the residual stress into the subassembly but also increase the cycling time (Yue et al. 2018, Du et al. 2019). Therefore, the optimal design of the fixture layout for a fuselage is important for reducing initial deformation, thus reducing the assembly cycling time and the residual stress.

Fixtures have been widely used in a variety of manufacturing systems to hold the parts to control their locations and orientations (Wang et al. 2010). Fixture design can be divided into several stages including setup planning, fixture planning, unit design, and verification (Wang et al. 2010). In this article, we focus on optimal fixture layout design in the fixture planning stage for complaint sheet parts.

In current industrial practice, Computer-aided fixture design (CAFD) methods can be categorized into three categories: rule-based methods, case-based reasoning (CRB) methods, and heuristic-based methods (Boyle et al. 2011). The first two categories belong to non-optimized layout planning methods (Boyle et al. 2011) which also highly rely on the re-use of experiential knowledge (Nee 1991, Kumar and Nee 1995, Joneja and Chang 1999, Zhang and Lin 1999, Gologlu 2004, Wang and Rong 2008). For the heuristic-based method, finite element analysis (FEA) tools are typically utilized to calculate the shape and deformation of the workpiece or subassembly with a given fixture layout. Then, a heuristic method is applied to search in the design space to find the fixture layout that minimizes the workpiece deformations. The commonly used heuristic methods include the Genetic Algorithm (GA) and Pseudo-gradient technique (Edward 1998, Li and Melkote 1999, Vallapuzha et al. 2002, Bazaraa et al. 2013). However, those methods all require a large amount of expensive FEA replications, which is time-consuming and even computationally prohibitive, especially in high precision modelling where a large number of mesh nodes is required. Recently, Du et al. (2021) proposed a method to tackle the optimal fixture layout design problem for compliant parts in the ship assembly. To avoid intensive implementations of FEA software, an integer programming problem was formulated using the stiffness matrix exported from the FEA software. It is a hard combinatorial optimization problem and the heuristic (simulated annealing (SA)) algorithm was adopted to solve it. In essence, this method is equivalent to running FEA software repeatedly and selecting the best result from multiple replications. As a result, the benefit of direct problem modeling using its stiffness matrix is not fully utilized and their approach still belongs to heuristic methods.

In this paper, we proposed a general fixture layout design framework (called SECR) which directly integrates system equations (SE) from the FEA software to formulate optimization problems, and the convex relaxation (CR) techniques are adopted to solve it. To the best knowledge of authors, this is the first non-heuristic-based general optimal fixture design framework.

There are two steps in the proposed SECR framework: In the first step, the system equation is exported from the FEA software, which describes the relationship between the linear shape deformation and forces applied on the part. This will serve as the foundation of our problem formulation. A detailed illustration of the system equation will be shown in section 2.1.

The second step is to formulate and solve the optimization problem. In a fixture design problem, finding optimal fixture locating point locations by using system equations is still intrinsically combinatorial. Convex relaxation is a powerful method in handling such problems. There are several famous convex relaxation techniques such as sparse learning relaxation (Donoho 2006, Du et al. 2019) and SDP relaxation (Goemans and Williamson 1995, (Sojoudi and Lavaei 2014, Beck et al. 2012) techniques. In particular, sparse learning relaxation is commonly used in transforming the combinatorial problem into a convex problem and has been proven to be effective in many engineering applications, such as compress sensing (Donoho 2006) and optimal actuator placement (Du et al. 2019). SDP relaxation is usually applied to tackle quadratic nonlinearities, which has been adopted in many areas such as graph theory (Goemans and Williamson 1995), approximation theory (Sojoudi and Lavaei 2014), and power systems (Beck et al. 2012). Compared to existing heuristic-based methods, our contribution is to derive a relaxation of the optimal fixture design problem by reformulating it into a convex SDP problem by using sparse learning relaxation and SDP relaxation techniques. The effectiveness of the convex relaxation algorithm is validated by a real case study in the half-to-half fuselage assembly process.

The remainder of this paper is organized as follows. Section 2 provides a detailed illustration of our convex relaxation-based optimal fixture design framework. A case study is provided to validate the performance of our proposed method in section 3. Finally, section 4 concludes the paper.

## 2. SECR Framework

In this section, we first elaborate on the FEA-based process model (Zhong et al. 2022). Then, we illustrate the general formulation for the optimal fixture design in section 2.2. Finally, the detailed derivation of convex relaxation is provided in section 2.3. In the convex relaxation, we utilize the sparsity property of fixture force which can further improve the computational efficiency and scalability. This makes it possible for general fixture design problems even in high precision modeling as the number of the nodes is not a concern in this method. By formulating the optimal fixture design problem in a convex manner, we can obtain an optimal solution by using the CVX software (Grant et al. 2008) efficiently.

## 2.1 FEA-based process model

For an assembly process that can be modelled in FEA software under static force and small linear deformation assumptions (Kohnke 2013), the relationship between linear elastic deformation and the force can be described by the following Equation (1):

$$\boldsymbol{K}\boldsymbol{U} = \boldsymbol{F}_g + \boldsymbol{F}_r, \tag{1}$$

where, $\boldsymbol{K} \in R^{6N\times 6N}$ is the global stiffness matrix; $\boldsymbol{U} = [\boldsymbol{u}_1; \dots; \boldsymbol{u}_N] \in R^{6N}$ is the nodal displacement vector assembled by aggerating the displacement vector $\boldsymbol{u}_i = \left[u_x^i, u_y^i, u_z^i, \omega_x^i, \omega_y^i, \omega_z^i\right]^T \in \mathbb{R}^{6\times 1}$ , $i \in \{1, .., N\}$, on each mesh node. Here, we denote the total number of mesh nodes as $N$ and a set containing all mesh nodes as $\mathcal{N}$. Moreover, we use $\boldsymbol{F}_r \in \mathbb{R}^{6N}$ and $\boldsymbol{F}_g \in \mathbb{R}^{6N}$ to represent the load exerted at the fixture locating points by the fixture and the gravity respectively. Vectors $\boldsymbol{F}_r$ and $\boldsymbol{F}_g$ are defined similarly as force vectors: $\boldsymbol{F}_r = [\boldsymbol{f}_{r1}; \dots; \boldsymbol{f}_{rN}] \in \mathbb{R}^{6N}$ and $\boldsymbol{F}_g = \left[\boldsymbol{f}_{g1}; \dots; \boldsymbol{f}_{gN}\right] \in \mathbb{R}^{6N}$ , where $\boldsymbol{f}_{gi} = \left[f_{gx}^i, f_{gy}^i, f_{gz}^i, \tau_{gx}^i, \tau_{gy}^i, \tau_{gz}^i\right]^T \in \mathbb{R}^6$ , $i \in \{1, .., N\}$ , is the gravity-induced load vector and $\boldsymbol{f}_{ri} = \left[f_{rx}^i, f_{ry}^i, f_{rz}^i, \tau_{rx}^i, \tau_{ry}^i, \tau_{rz}^i\right]^T \in \mathbb{R}^6$ , $i \in \{1, .., N\}$, is the fixture-induced load vector on each mesh node. Specifically, the first three elements of $\boldsymbol{f}_{gi}$, $\left[f_{gx}^i, f_{gy}^i, f_{gz}^i\right]^T$, denote the three-dimensional load-induced force of the $i$th mesh node; and the last three elements of $\boldsymbol{f}_{gi}$, $\left[\tau_{gx}^i, \tau_{gy}^i, \tau_{gz}^i\right]^T$, denote the three-dimensional load-induced torque of the $i$th mesh node. $\boldsymbol{f}_{ri}$ is defined similarly as $\boldsymbol{f}_{gi}$ for its elements. Without loss of generality, we assume fixtures can only constraint linear displacement, i.e., $\left[u_x^i, u_y^i, u_z^i\right]^T = \boldsymbol{0}$, when the $i$th mesh node is used as a fixture locating point. If a fixture can constraint both linear and angular displacement, then, $\left[u_x^i, u_y^i, u_z^i, \omega_x^i, \omega_y^i, \omega_z^i\right]^T = \boldsymbol{0}$, when the $i$th mesh node is used as a fixture locating point.

The stiffness matrix $\boldsymbol{K}$ and force $\boldsymbol{F}_g$ can be exported from the FEA software directly. To apply our algorithm, we assume that at least three locating points that make the structure stable are pre-specified. Once those three fixtures locating points are fixed, the fuselage has unique deformation under gravity. Then, we follow the common practice in FEA (Kohnke 2013) to remove the corresponding rows and columns of those three pre-specified fixture mesh nodes in $\boldsymbol{K}$ to obtain an invertible stiffness matrix $\boldsymbol{K}^* \in \mathbb{R}^{6N_1 \times 6N_1}$, where $N_1 = N - 3$. Denote the remaining set of nodes as $\mathcal{N}_1$. Similarly, by removing the corresponding rows, we can obtain $\boldsymbol{U}^*$, $\boldsymbol{F}_g^*$ and $\boldsymbol{F}_r^*$, respectively. Then, by conducting basic linear algebra operations, we obtain:

$$\boldsymbol{K}^*\boldsymbol{U}^* = \boldsymbol{F}_g^* + \boldsymbol{F}_r^*.$$

Since $\boldsymbol{K}^*$ is an invertible matrix, we use $\boldsymbol{A}^*$ to denote the inverse matrix of $\boldsymbol{K}^*$:

$$\boldsymbol{U}^* = \boldsymbol{A}^*\left(\boldsymbol{F}_g^* + \boldsymbol{F}_r^*\right) \tag{2}$$

## 2.2 Original formulation for optimal fixture design

For a complaint part, we adopt 'N-2-1' locating principle in the fixture design (Cai et al. 1996). Denote the set of potential fixture locating points as $\mathcal{N}_{PT} \subseteq \mathcal{N}_1$ and its cardinality as $N_{PT}$. We aim to find the set of $N_f$ fixture locating points, denoted as $\mathcal{N}_f$, out of $N_{PT}$ potential fixture locating points to achieve minimum total deformation $\delta^2$, i. e.,

$$\delta^2 = (\boldsymbol{U}^*)^T \boldsymbol{W} \boldsymbol{U}^*,$$

where $\boldsymbol{U}^*$ is the shape deformation vector; $\boldsymbol{W}$ is a diagonal matrix with only ones on locations of the linear displacement $\left[u_x^i, u_y^i, u_z^i\right]^T$ of mesh nodes that we are interested in, and all zeros otherwise, i.e.,

$$diag(\boldsymbol{W}) = \left[0,0,0,\dots,\underbrace{1,1,1,0,0,0}_{i}\dots,0,0,0\right]^T,$$

where $diag(\boldsymbol{W})$ returns the diagonal elements of matrix $\boldsymbol{W}$.

In reality, the following physical constraints are considered in the fixture design optimization:

- Notice that compared to the number of nodes, the number of fixture locating points is usually much smaller. Therefore, $\boldsymbol{F}_r^*$ is a sparse vector whose nonzero elements can only appear in potential fixture locating points. Let $\overline{\mathcal{N}}_{PT}$ represent the complementary set of $\mathcal{N}_{PT}$, i.e., $\overline{\mathcal{N}}_{PT} = \mathcal{N}_1 \backslash \mathcal{N}_{PT}$. Let the fixture induced force vector at node $i$ be $\boldsymbol{f}_i^* = \left[f_{rx}^i, f_{ry}^i, f_{rz}^i\right]^T$, i.e., $\boldsymbol{f}_i^* = \boldsymbol{0}, \forall i \in \overline{\mathcal{N}}_{PT}$. Let $\boldsymbol{F}_{rPT}^* = \left[\boldsymbol{f}_{ri_1}^*; \dots; \boldsymbol{f}_{ri_l}^*; \dots; \boldsymbol{f}_{ri_{PT}}^*\right]^T$, where $i_l \in \mathcal{N}_{PT}$, $l \in \{1, \dots, N_{PT}\}$. Notice that each element in $\boldsymbol{F}_{rPT}^*$ has its corresponding index in $\boldsymbol{F}_r^*$. In order to link the indices of elements between $\boldsymbol{F}_{rPT}^*$ and $\boldsymbol{F}_r^*$, we define the following mapping: the $l$-th element in $\boldsymbol{F}_{rPT}^*$ is the $\left(i_{(floor(l/3)+1)} - 1\right) \times 6 + (l \, mod \, 3)$-th element in $\boldsymbol{F}_r^*$, where $i_k \in \mathcal{N}_{PT}$, $k \in \{1, \dots, N_{PT}\}$, $floor(l/3)$ denotes the quotient and $l \, mod \, 3$ denotes the remainder. Define this mapping as $m(l) = \left(i_{(floor(l/3)+1)} - 1\right) \times 6 + (l \, mod \, 3)$.
- Let $\boldsymbol{\tau}_i^*$ be the fixture induced torque vector on mesh node $i$. Note that fixtures only exert force and no torque on the part, $\boldsymbol{\tau}_i^*$ on all mesh nodes should be zero, i.e., $\boldsymbol{\tau}_i^* = \left[\tau_{rx}^i, \tau_{ry}^i, \tau_{rz}^i\right]^T = \boldsymbol{0}, i \in \mathcal{N}_1$.
- Note that the fixture locating point will restrict the deviation of the corresponding mesh node to be zero, i.e., $\left[u_x^i, u_y^i, u_z^i\right]^T = \boldsymbol{0}, i \in \mathcal{N}_f$. For non-fixture location points, $\left[f_{rx}^i, f_{ry}^i, f_{rz}^i\right]^T = \boldsymbol{0}$, $i \in \overline{\mathcal{N}_f}$. Here, $\overline{\mathcal{N}_f}$ represents the non-fixture location points. Then, for any given points $i$, we will always have $u_a^i f_{ra}^i = 0, i \in \{1, \dots, N_1\}, a \in \{x, y, z\}$. Moreover, since $\boldsymbol{\tau}_i = \boldsymbol{0}, i \in \mathcal{N}_1$, we can further write this constraint as $\boldsymbol{U}^*(l)\boldsymbol{F}_r^*(l) = 0, j \in \{1, \dots, 6N_1\}$. If we focus on potential fixture locations, we have $\boldsymbol{U}^*\big(m(l)\big)\boldsymbol{F}_{rPT}^*(l) = 0, l \in \{1, \dots, 3N_{PT}\}$.
- We further define the magnitudes of fixture induced force vector on all potential fixture locations as $\boldsymbol{F}_1^* = \left[\|\boldsymbol{f}_1^*\|_2, \dots, \left\|\boldsymbol{f}_{N_1}^*\right\|_2\right]^T$. When we have $n_a$ fixture locating points, there are only $n_a$ fixture force vectors with nonzero elements. Thus, we have $\|\boldsymbol{F}_1\|_0 = n_a$, where $\|\cdot\|_0$ represents the number of nonzero entries in a vector.

Finally, the problem can be formulated as

$$\min_{\boldsymbol{F}_{rPT}^*} (\boldsymbol{U}^*)^T \boldsymbol{W} \boldsymbol{U}^* \tag{3}$$

$$\text{subject to: } \boldsymbol{U}^* = \boldsymbol{A}^*\left(\boldsymbol{F}_g^* + \boldsymbol{F}_r^*\right) \tag{3a}$$

$$\boldsymbol{f}_i^* = \boldsymbol{0}, i \in \overline{\mathcal{N}}_{PT} \tag{3b}$$

$$\boldsymbol{\tau}_i^* = \boldsymbol{0}, i \in \mathcal{N}_1 \tag{3c}$$

$$\left[u_x^i, u_y^i, u_z^i\right]^T = \boldsymbol{0}, i \in \mathcal{N}_f \tag{3d}$$

$$\boldsymbol{\tau}_i^* = \left[\tau_{rx}^i, \tau_{ry}^i, \tau_{rz}^i\right]^T, \boldsymbol{f}_i^* = \left[f_{rx}^i, f_{ry}^i, f_{rz}^i\right]^T, \boldsymbol{f}_{ri}^* = [\boldsymbol{\tau}_i^*;\ \boldsymbol{f}_i^*]\ i \in \mathcal{N}_1 \text{ and } \boldsymbol{F}_r^* = \left[\boldsymbol{f}_{r1}^*; \dots; \boldsymbol{f}_{rN_1}^*\right] \tag{3e}$$

$$\boldsymbol{F}_1^* = \left[\|\boldsymbol{f}_1^*\|_2, \dots, \left\|\boldsymbol{f}_{N_1}^*\right\|_2\right]^T \tag{3f}$$

$$\|\boldsymbol{F}_1^*\|_0 = n_a \tag{3g}$$

$$\boldsymbol{U}^*\big(m(l)\big)\boldsymbol{F}_{rPT}^*(l) = 0, l \in \{1, \dots, 3N_{PT}\} \tag{3h}$$

$$\boldsymbol{F}_{rPT}^* = \left[\boldsymbol{f}_{ri_1}^*; \dots; \boldsymbol{f}_{ri_l}^*; \dots; \boldsymbol{f}_{ri_{PT}}^*\right]^T, i_l \in \mathcal{N}_{PT}, l \in \{1, \dots, N_{PT}\} \tag{3$i$}$$

## 2.3 Convex formulation for optimal fixture design

The selection of the set of $\mathcal{N}_f$ fixture locating points out of mesh nodes can be formulated as an integer programming problem. Previous research demonstrated that traditional integer programming methods are hard to be applied (Du et al. 2021) in fixture design problems due to the huge number of mesh nodes or too large search space. Therefore, heuristic methods are prevalent in this field. To avoid the drawback of heuristic methods, we propose a convex relaxation method for solving the problem (3). The main challenge is to find a good convex relaxation of the following non-convex constraints: (a) the $l_0$ norm constraint (3g), and (b) complementary slackness constraint (3h).

*(a) Convex relaxation of $l_0$ norm*

To handle the $l_0$ norm constraint (3g), we follow the similar procedure of (Du et al. 2019) and transform the optimization problem (3) into an optimization problem with group lasso penalty:

$$\min_{\boldsymbol{F}_{rPT}^*} (\boldsymbol{U}^*)^T \boldsymbol{W} \boldsymbol{U}^* + \lambda \|\boldsymbol{F}_1^*\|_1 \tag{4}$$

$$\text{subject to: Constraints (3a)} - \text{(3e), (3h), (3}i\text{)}$$

*(b) Convex relaxation of complementary slackness constraint*

Finding a good convex relaxation of the complementary slackness constraint (3h) is tricky. There are mainly three steps in relaxing this constraint: (i) we first exploit the linear relationship between $\boldsymbol{U}^*$ and $\boldsymbol{F}_r^*$ to transform the complementary slackness constraint into a quadratic equality constraint; (ii) then, we further transform this problem into an SDP problem with rank-1 constraint; (iii) finally, we relax this rank-1 constraint by penalizing its nuclear norm and thereby transform the problem (4) into a convex problem. The relaxation steps are presented as follows:

<u>(i) Transforming constraint (3h) into a quadratic equality constraint.</u> Since the linear relationship between $\boldsymbol{U}^*$ and $\boldsymbol{F}_r^*$ is known in constraint (3a), the quadratic equality constraint can be obtained by substituting $\boldsymbol{U}^*$ with $\boldsymbol{F}_r^*$ . Notice that, either substituting $\boldsymbol{U}^*$ with $\boldsymbol{F}_r^*$ or substituting $\boldsymbol{F}_r^*$ with $\boldsymbol{U}^*$ are mathematically equivalent. However, we choose the first option to utilize the sparsity property of $\boldsymbol{F}_r^*$, which highly improves the computational efficiency and makes the method scalable. Adopting this approach, the scale of the SDP problem is irrelevant to the number of nodes which can be seen in the following derivations. After substituting $\boldsymbol{F}_r^*$ with $\boldsymbol{U}^*$, (3h) can be written as

$$\sum_{j=1}^{6N_1} \boldsymbol{A}^*(m(l),j)\boldsymbol{F}_g^*(j)\boldsymbol{F}_{rPT}^*(l) + \sum_{k\in\{1,\dots,3N_{PT}\}} \boldsymbol{A}^*\big(m(l),m(k)\big)\boldsymbol{F}_r^*\big(m(k)\big)\boldsymbol{F}_r^*\big(m(l)\big) = 0, l \in \{1,\dots,3N_{PT}\}. \quad (3h')$$

<u>(ii) Transforming quadratic equality constraint into SDP problem with a rank-1 constraint:</u> The constraint (3h′) is still nonconvex. To deal with this quadratic equality constraint, we define a new matrix $\boldsymbol{S_1} = \boldsymbol{F}_{rPT}^*(\boldsymbol{F}_{rPT}^*)^T$. Then we have $\boldsymbol{S}_1(l,k) = \boldsymbol{F}_{rPT}^*(l)\boldsymbol{F}_{rPT}^*(k) = \boldsymbol{F}_r^*\big(m(l)\big)\boldsymbol{F}_r^*\big(m(k)\big)$ and constraint (3h′) is equivalent to

$$\sum_{j=1}^{6N_1} \boldsymbol{A}^*(m(l),j)\boldsymbol{F}_g^*(j)\boldsymbol{F}_{rPT}^*(l) + \sum_{k\in\{1,\dots,3N_{PT}\}} \boldsymbol{A}^*\big(m(l),m(k)\big)\boldsymbol{S}_1(l,k) = 0, l \in \{1,\dots,3N_{PT}\}. \quad (3h'')$$

Now, constraint (3h″) is a linear constraint. However, constraint $\boldsymbol{S_1} = \boldsymbol{F}_{rPT}^*(\boldsymbol{F}_{rPT}^*)^T$ is still nonconvex, which makes solving problem (4) intractable.

To solve this problem, an idea is to reformulate problem (4) into an SDP problem. Next, the following theorem establishes the equivalence between the original problem (4) and an SDP problem with rank-1 constraint.

**Theorem 1.** If there is an $l \in \{1, \dots, 3N_{PT}\}$ such that $\boldsymbol{F}_{rPT}^*(l)\sum_{j=1}^{6N_1}\boldsymbol{A}^*(m(l),j)\boldsymbol{F}_g^*(j) \neq 0$, then problem (4) is equivalent to the problem (5) below:

$$\min_{\boldsymbol{F}_{rPT}^*,\ \boldsymbol{S}_1} (\boldsymbol{U}^*)^T\boldsymbol{W}\boldsymbol{U}^* + \lambda\|\boldsymbol{F}_1^*\|_1 \tag{5}$$

$$\text{subject to: Constraints (3a)} - \text{(3e), (3i), (3h}''\text{)}$$

$$\boldsymbol{S}_1 - \boldsymbol{F}_{rPT}^*(\boldsymbol{F}_{rPT}^*)^T \succcurlyeq \boldsymbol{0} \tag{5o}$$

$$\text{rank}(\boldsymbol{S}_1) = 1 \tag{5r}$$

The proof of Theorem 1 is given in the Appendix. The assumption in Theorem 1 is easy to check. We propose two approaches: (i) before solving the problem: to balance the gravity load, at least one element in potential fixture locations is non-zero (the location is unknown). In this case, if $\sum_{j=1}^{6N_1}\boldsymbol{A}^*(m(l),j)\boldsymbol{F}_g^*(j) \neq 0, \forall\, l \in \{1, \dots, 3N_{PT}\}$, then the assumption is valid. (ii) after solving the problem: the locations ($l$) of non-zero elements in potential fixture locations are known. In this case, if at least one $l$ such that $\boldsymbol{F}_{rPT}^*(l)\sum_{j=1}^{6N_1}\boldsymbol{A}^*(m(l),j)\boldsymbol{F}_g^*(j) \neq 0$, then the assumption is valid.

Notice that the dimension of this SDP problem is only related to the number of potential fixtures locating points, which is usually small and independent of the number of mesh nodes. There are some existing convex relaxation methods proposed to deal with constraint (5r), i.e., nuclear norm relaxation method (Zhang et al. 2012, Mou et al. 2020). Before dealing with constraint (5r), we can write constraint (5o) more compactly. To achieve this, we define $\boldsymbol{S}_2 = \begin{bmatrix} \boldsymbol{S}_1 & \boldsymbol{F}_{rPT}^* \\ (\boldsymbol{F}_{rPT}^*)^T & 1 \end{bmatrix}$ and constraint (5o) is transformed into:

$$\mathbf{S}_2 \succcurlyeq \mathbf{0} \tag{5$o'$}$$

$$\boldsymbol{S}_2 = \begin{bmatrix} \boldsymbol{S}_1 & \boldsymbol{F}_{rPT}^* \\ (\boldsymbol{F}_{rPT}^*)^T & 1 \end{bmatrix} \tag{5$o''$}$$

(iii) SDP relaxation by applying nuclear norm penalty: By applying the nuclear norm relaxation method, we can remove the constraint (5r) by adding the nuclear norm penalty in the objective function. Finally, we transform problem (5) into the problem (6) below:

$$\min_{\boldsymbol{F}_{rPT}^*,\ \boldsymbol{S}_2} (\boldsymbol{U}^*)^T \boldsymbol{W} \boldsymbol{U}^* + \lambda \|\boldsymbol{F}_1^*\|_1 + \mu \|\boldsymbol{S}_2\|_* \tag{6}$$

$$\text{subject to: } \boldsymbol{U}^* = \boldsymbol{A}^*\left(\boldsymbol{F}_g^* + \boldsymbol{F}_r^*\right) \tag{6a}$$

$$\boldsymbol{f}_i^* = \boldsymbol{0}, i \in \bar{\mathcal{N}}_{PT} \tag{6b}$$

$$\boldsymbol{\tau}_i^* = \boldsymbol{0}, i \in \mathcal{N}_1 \tag{6c}$$

$$\left[u_x^i, u_y^i, u_z^i\right] = \boldsymbol{0}, i \in \mathcal{N}_f \tag{6d}$$

$$\boldsymbol{\tau}_i^* = \left[\tau_{rx}^i, \tau_{ry}^i, \tau_{rz}^i\right], \boldsymbol{f}_i^* = \left[f_{rx}^i, f_{ry}^i, f_{rz}^i\right]^T, \boldsymbol{f}_{ri}^* = \left[\boldsymbol{\tau}_i^*;\ \boldsymbol{f}_i^*\right]\ i \in \mathcal{N}_1 \text{ and } \boldsymbol{F}_r^* = \left[\boldsymbol{f}_{r1}^*; \dots; \boldsymbol{f}_{rN_1}^*\right] \tag{6e}$$

$$\boldsymbol{F}_1^* = \left[\|\boldsymbol{f}_1^*\|_2, \dots, \left\|\boldsymbol{f}_{N_1}^*\right\|_2\right]^T \tag{6f}$$

$$\sum_{j=1}^{6N_1} \boldsymbol{A}^*(m(l), j) \boldsymbol{F}_g^*(j) \boldsymbol{F}_{rPT}^*(l) + \sum_{l \in \{1,\dots,3N_{PT}\}} \boldsymbol{A}^*\big(m(l), m(k)\big) \boldsymbol{S}_1(l,k) \ = 0, l \in \{1, \dots, 3N_{PT}\} \tag{6h}$$

$$\boldsymbol{F}_{rPT}^* = \left[\boldsymbol{f}_{ri_1}^*; \dots; \boldsymbol{f}_{ri_l}^*; \dots; \boldsymbol{f}_{ri_{PT}}^*\right], i_l \in \mathcal{N}_{PT}, l \in \{1, \dots, N_{PT}\} \tag{6l}$$

$$\boldsymbol{S}_2 \succcurlyeq \boldsymbol{0} \tag{6o}$$

$$\boldsymbol{S}_2 = \begin{bmatrix} \boldsymbol{S}_1 & \boldsymbol{F}_{rPT}^* \\ (\boldsymbol{F}_{rPT}^*)^T & 1 \end{bmatrix} \tag{6p}$$

Problem (6) is the final problem formulation, which is a typical SDP problem and can be solved efficiently by using the CVX software. We find $\lambda$ value by using a binary search algorithm proposed in (Du et al. 2019). Here is the general idea: if the tuning parameter $\lambda$ is too large, we will have $\|\boldsymbol{F_1}\|_0 < n_a$, then we should decrease the $\lambda$ value; when the tuning parameter $\lambda$ is too small, we will have $\|\boldsymbol{F_1}\|_0 > n_a$, then we should increase the $\lambda$ value. The selection of the tuning parameter $\mu$ is based on the rank of $\boldsymbol{S_2}$ matrix, when the rank of $\boldsymbol{S_2}$ is not equal to 1, we should increase the value of $\mu$ until that the constraint $rank(\boldsymbol{S_2}) = 1$ is satisfied.

## 3. Case study

In this case study, we will use a half-to-half fuselage assembly process to demonstrate the proposed SECR framework. In the fuselage assembly process, gravity-induced deformation is directly related to assembly precision. Different fixture layouts will result in a maximum shape deformation induced by gravity load varying from less than 0.01 inches to larger than 10 inches. For example, we have two different sets of fixtures in Figures 1 (a) and (b). For the fixture layout shown in Figure 1 (a), the maximum total deformation is 0.02 inches, and the maximum residual stress is 611.3 psi. While for fixture layout in Figure 1(b), the maximum total deformation is 0.13 inches, and the maximum residual stress is 1497.3 psi. From Figure 1, we can see that the fixture layout will significantly influence the precision of the half fuselage assembly process. However, both layouts cannot meet the engineering specifications since ultra-high precision assembly is vital for large-scale aircraft production. Without a well-designed fixture layout, it is challenging to achieve such high precision in assembly. A proper fixture layout should be able to compensate for the gravity-induced shape deformation. In this application, our objective is to find the optimal locations for a given number of fixture locating points such that the gravity-induced deformation is minimized.

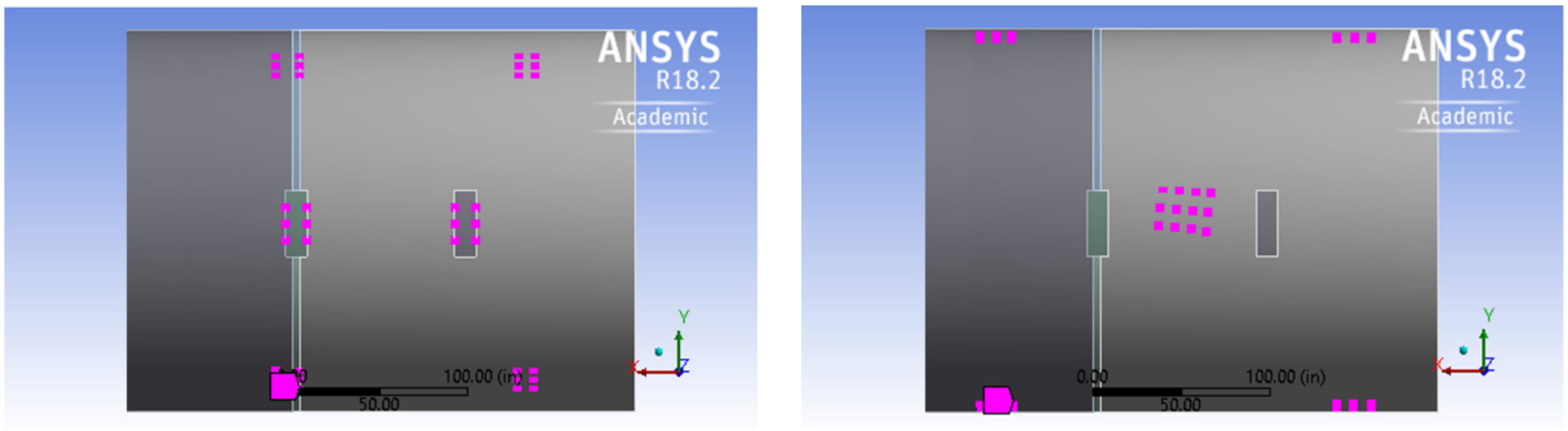

(a) Fixture layout 1 (b) Fixture layout 2

Figure 1. Two different fixture allocation schemes.

## 3.1 Comparison study of the optimally designed fixture with the current industrial practice

In the proposed SECR method, we assume that there are 30 locations for the potential fixture locating points and 3 pre-specified locations for the fixture locating points, which are shown in Figure 2. In Figure 2, the red points and the green points represent the mesh nodes and locations for potential fixture locating points, respectively. We checked that the assumptions in Theorem 1 hold. The tuning parameters are selected as $\lambda = 5$ and $\mu = 10$.

To show the superiority of the proposed method, we first conduct a comparison among the results of the proposed method and the current industry practice. The fixture layout corresponding to the proposed method is shown in Figure 3(a), while the fixture layout of current industry practice is shown in Figure 3(b). With the same number of fixture locating points (e.g. both have $N_f = 8$), the SECR method can achieve maximum total deformation of 0.014 inches while the maximum total deformation of current industry practice is 0.27 inches. This result indicates that with the same number of fixture locating points, our algorithm can significantly reduce the maximum total deformation by optimizing the fixture locating point locations.

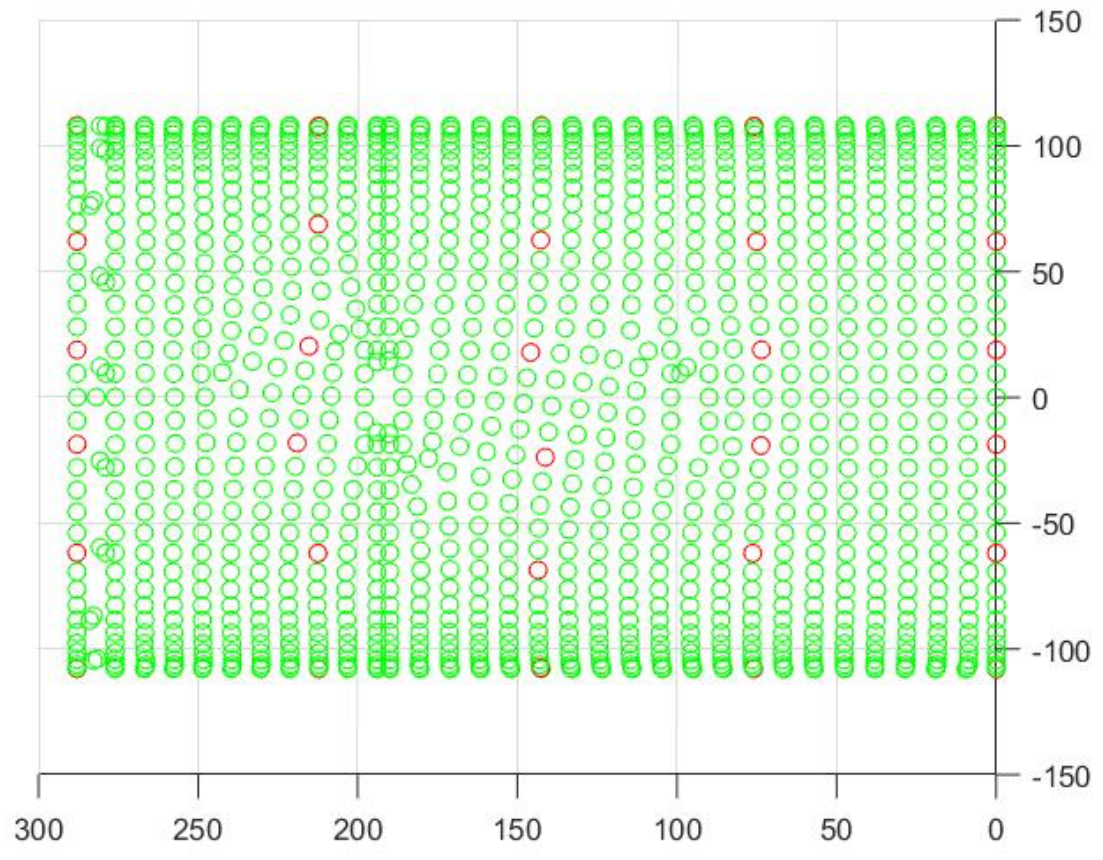

Figure 2. Potential locations to place fixture locators: red dot –potential locations; green dot – mesh nodes

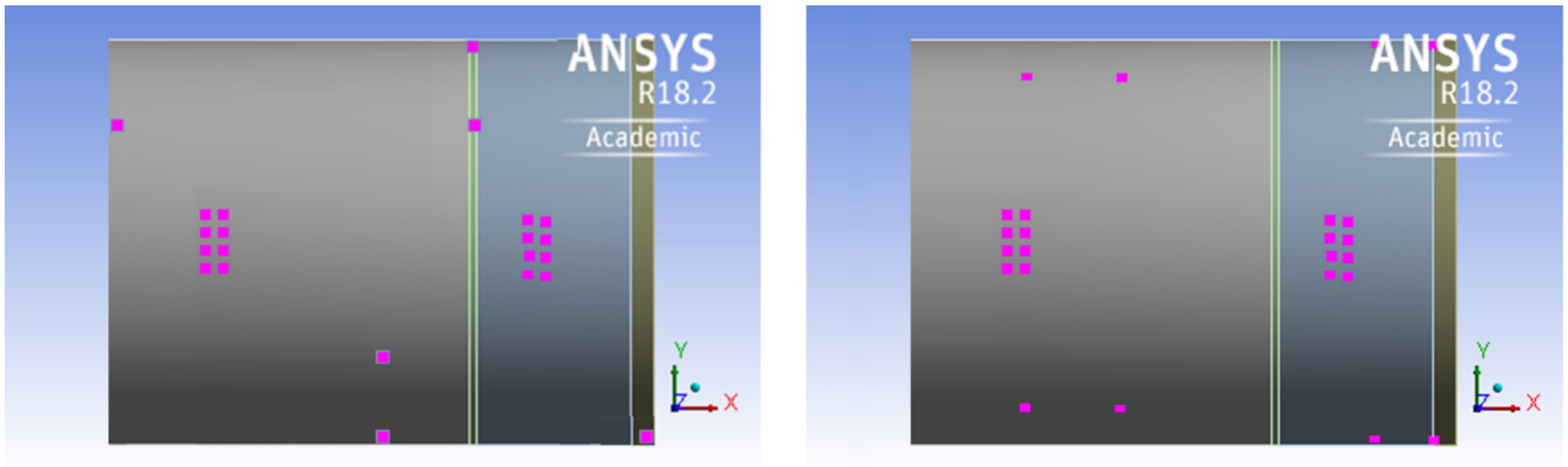

(a) Fixture layouts of the SECR algorithm (b) Fixture layouts of current industry practice

Figure 3. Fixture layouts of the SECR algorithm and current industry practice.

### *3.2 Comparison of the optimal design with heuristic-based design*

We also compare the result of the proposed method with the heuristic-based DSMSO method (Du et al. 2021). The method in Du et al. (2021) first loads the stiffness matrix from the FEA simulation platform and an integer programming problem aiming at minimizing the maximum deformation is formulated. Then, simulated annealing is adopted to solve this integer programming problem.

We conduct the following two comparisons:

1. We run DSMSO methods for 300 replications with the same amount of computational time as our methods. By using the same settings as above, the maximum total deformation of our fixture layouts is 0.0112 inches. The best scenario maximum total deformation of the DSMSO method among those 300 replications is 0.0159 inches. The comparison results are documented in Figure 4. In Figure 4, the x-axis represents the replications while the y-axis represents the maximum total deformation. The curve above represents the result of the DSMSO method, and the curve below stands for the result of the proposed SECR method. The result indicates that the proposed method outperforms the state-of-art DSMSO method.

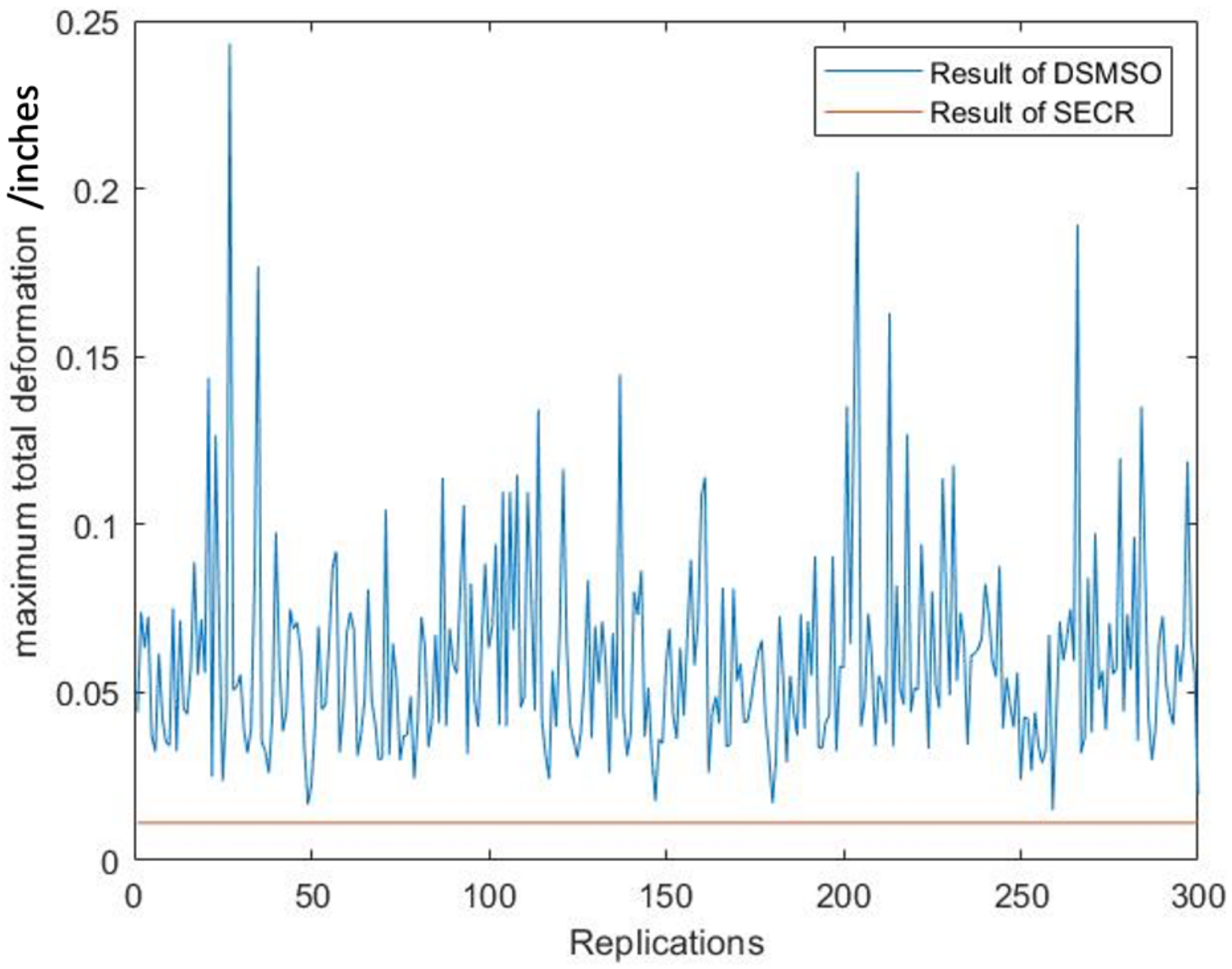

Figure 4. Comparison results with DSMSO method.

2. We run the DSMSO methods for 1000 replications with the same SA stopping criteria proposed by Du et al. (2021). As a result, each SA replication has 1044 FEA runs and takes 2016 seconds on average, which is approximately 30 times of the SECR running time. Figure 5 shows the histogram of the maximum total deformation of those 1000 DSMSO replications.

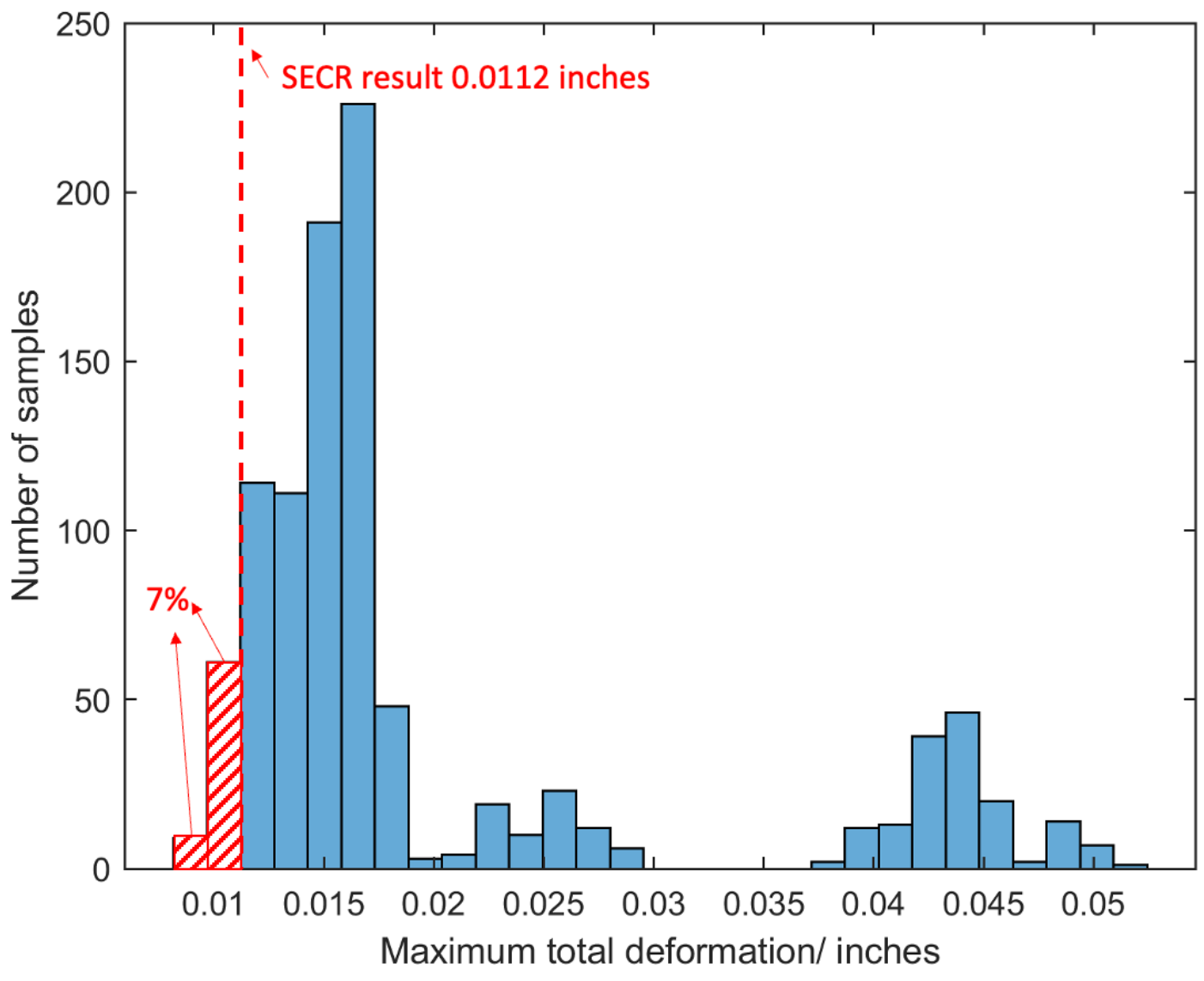

Figure 5. Histogram of the maximum total deformation of DSMSO replications.

The x-axis represents the maximum total deformation and y-axis represents the number of DSMSO replication results in a specific maximum total deformation range. It approximates the distribution of the DSMSO result. Among those 1000 DSMSO replications, only 7% achieve the same (or slightly better) result than that of the SECR method. This means, on average, 14.3 DSMSO replications (8 hours) are needed to achieve a similar result as the SECR method which only takes 69 seconds. Moreover, the SECR result can serve as a warm start for the DSMSO method, which will significantly reduce the computation time and improve the solution quality.

We also tried to use the solver Gurobi to solve the problem (3). However, the solver does not make any progress in 10 hours. Therefore, we terminated it.

Another advantage of the proposed method is its computational efficiency. For heuristic-based methods, a large linear system has to be solved repeatedly and there is no guarantee of the optimality of its solution. However, by using our proposed methods, we can obtain the results by simply solving a small-scale SDP problem, which is much more computationally efficient. For comparison, we also document the computational time for a different number of potential fixture locations in Table 1. As can be seen from the table, the computational time of solving the SDP problem will increase significantly as the scale (number of potential fixture locations) increases, since the complexity of solving such an SDP problem per iteration is $O(N_{PT}^3)$. This shows the necessity of adopting sparsity property to reduce the problem scale in the proposed SECR method. In real practice, by simply choosing 30-40 potential fixture locations, we can already obtain a result within tolerance.

**Table 1. The computational time of SECR for different number of potential fixture locations**

| Number of potential fixture locations | Computational time |
|---|---|
| 20 | 29s |
| 30 | 69s |
| 42 | 358s |
| 56 | 1661s |

## 4. Conclusion

This paper proposed a novel SECR framework for general fixture layout design. Since the optimal fixture layout design is generally a large-scale combinatorial optimization problem, convex-relaxation techniques, especially sparse-learning and SDP relaxation techniques, are adopted to formulate the optimal fixture design problem in a convex manner. The proposed framework is computationally efficient and scalable. The solution of this convex problem implies a near-optimal fixture layout, and this convex problem can be solved efficiently by using existing convex optimization algorithms such as interior-point methods, ADMM, etc.

To validate the effectiveness of the proposed method, we conducted a real case study in a half-to-half fuselage assembly process. In the case study, we compare the fixture layout generated by our algorithm with the current industrial practice and the DSMSO method. The result shows that the fixture layout obtained by the proposed SECR method outperforms both the current fixture used in a real assembly process and the fixture layout generated by the DSMSO method in terms of maximum total deformation.

Our proposed SECR algorithm was developed with the assumptions of small linear deformations and static forces, which can be applied to broad applications beyond fuselage fixture layout design. Because

- In many applications, the deformation is restricted to the linear elastic region to avoid material failure;
- The static force can be static manufacturing or assembly forces instead of gravity.

The SECR framework can also be generalized to broad settings by considering either relaxing the small linear deformation assumption or static force assumptions.

## Data availability statement

Due to the nature of this research, participants of this study did not agree for their data to be shared publicly, so supporting data is not available.

## Acknowledgment

The work is supported by the Strategic University Partnership between the Boeing Company and the Georgia Institute of Technology (Funder ID: 10.13039/ 100000003).

## Nomenclature

$\mathcal{N}$ = a set contains all mesh node

N = size of set $\mathcal{N}$

$\boldsymbol{K}$ = global stiffness matrix

$u_k^i$ = the linear displacement of the $i$th mesh node in direction $k$, where $k = x, y, z$.

$\omega_k^i$ = the angular displacement of the $i$ th mesh node in direction $k$, where $k = x, y, z$.

$\boldsymbol{u}_i$ = nodal displacement vector on the $i$ th mesh node defined as $\left[u_x^i, u_y^i, u_z^i, \omega_x^i, \omega_y^i, \omega_z^i\right]^T, i \in \{1, .., N\}$.

$\boldsymbol{U}$ = nodal displacement vector of all mesh nodes defined as $[\boldsymbol{u}_1; \dots; \boldsymbol{u}_N]$

$f_{gk}^i$ = gravity-induced force on the $i$ th mesh node in direction $k$, where $k = x, y, z$.

$\tau_{gk}^i$ = gravity-induced torque on the $i$ th mesh node in direction $k$, where $k = x, y, z$.

$\boldsymbol{f}_{gi}$ = gravity-induced load vector on the $i$th mesh node defined as $\boldsymbol{f}_{gi} = \left[f_{gx}^i, f_{gy}^i, f_{gz}^i, \tau_{gx}^i, \tau_{gx}^i, \tau_{gz}^i\right]^T, i \in \{1, .., N\}$

$\boldsymbol{F_g}$ = gravity-induced load vector of all mesh nodes defined as $\left[\boldsymbol{f}_{g1}; \dots; \boldsymbol{f}_{gN}\right]$.

$f_{rk}^i$ = fixture locating points induced force on the $i$ th mesh node in direction $k$, where $k = x, y, z$.

$\tau_{rk}^i$ = fixture locating points induced torque on the $i$ th mesh node in direction $k$, where $k = x, y, z$.

$\boldsymbol{f}_{ri}$ = gravity-induced load vector on the $i$ th mesh node defined as $\boldsymbol{f}_{ri} = \left[f_{rx}^i, f_{ry}^i, f_{rz}^i, \tau_{rx}^i, \tau_{ry}^i, \tau_{rz}^i\right]^T, i \in \{1, .., N\}$.

$\boldsymbol{F_r}$ = fixture locating points induced load on all mesh nodes defined as $[\boldsymbol{f}_{r1}; \dots; \boldsymbol{f}_{rN}]$

$\mathcal{N}_1$ = a set contain all mesh nodes beside three fixture locating points

$N_1$ = size of set $\mathcal{N}_1$

$\boldsymbol{K}^*$ = stiffness matrix removing the corresponding rows and columns of fixture locating points

$\boldsymbol{A}^* = (\boldsymbol{K}^*)^{-1}$

$\boldsymbol{U}^*$ = nodal displacement vector after removing the corresponding rows of prespecified fixture locating points

$\boldsymbol{F}_g^*$ = gravity-induced load vector of all mesh nodes after removing the corresponding rows of pre-specified fixture locating points

$\boldsymbol{f}_i^*$ = fixture induced 3-dimensional force vector on the $i$th mesh node defined as $\left[f_{rx}^1, f_{ry}^1, f_{rz}^1\right]^T, i \in \{1, \dots, N^*\}$

$\boldsymbol{\tau}_i^*$ = fixture induced 3-dimensional torque vector on the $i$th mesh node defined as $\left[\tau_{rx}^i, \tau_{ry}^i, \tau_{rz}^i\right]^T, i \in \{1, \dots, N^*\}$

$\boldsymbol{f}_{ri}^*$ = fixture induced load vector on the $i$th mesh node defined as $\boldsymbol{f}_{ri}^* = [\boldsymbol{f}_i^*; \boldsymbol{\tau}_i^*], i \in \{1, \dots, N^*\}$

$\boldsymbol{F}_r^*$ = fixture induced load vector of all mesh nodes after removing the corresponding rows of prespecified fixture locating points, defined as $[\boldsymbol{f}_{r1}^*; \dots; \boldsymbol{f}_{rN^*}^*]$

$\boldsymbol{W}$ = diagonal matrix with ones and all other elements are zeros

$\delta^2$ = total deviation

$\mathcal{N}_{\mathrm{PT}}$ = a set of potential fixture locations

$N_{PT}$ = size of $\mathcal{N}_{\mathrm{PT}}$

$\bar{\mathcal{N}}_{PT}$ = complementary set of potential fixture locating points locations

$\mathcal{N}_f$ = the set of fixture locating points locations

$\overline{\mathcal{N}_f}$ = represents the non-fixture location points

$\boldsymbol{F}_1^*$ = the magnitudes of fixture induced force vector defined as $\left[\|\boldsymbol{f}_1^*\|_2, \dots, \left\|\boldsymbol{f}_{N_1}^*\right\|_2\right]^T$

$n_a$ = number of available fixtures locating points

$\boldsymbol{S_0} = \boldsymbol{F}_r^*(\boldsymbol{F}_r^*)^T$

$\boldsymbol{F}_{rPT}^* = [\boldsymbol{f}_{ri_1}^*, \dots, \boldsymbol{f}_{ri_l}^*, \dots, \boldsymbol{f}_{ri_{PT}}^*]$

$\boldsymbol{S_1} = \boldsymbol{F}_{rPT}^*(\boldsymbol{F}_{rPT}^*)^T$

$\boldsymbol{S_2} = \begin{bmatrix} \boldsymbol{S}_1 & \boldsymbol{F}_{rPT}^* \\ (\boldsymbol{F}_{rPT}^*)^T & 1 \end{bmatrix}$

$\lambda, \mu$ = tuning parameters

## Appendix. Proof of Theorem 1.

Proof: Since $(5)$ is the convex relaxation of $(4)$, the solution of $(4)$ is also a feasible solution of $(5)$. We will now prove that $(5)$ implies $(4)$:

(i) Since $rank(\boldsymbol{S}_1) = 1$, $\boldsymbol{S}_1$ can be written as $\boldsymbol{S_1} = \boldsymbol{v}\boldsymbol{v}^{\boldsymbol{T}}$, where $\boldsymbol{v} \in \mathbb{R}^{3N_{PT}}$. Constraint (5o) can be written as $\boldsymbol{v}\boldsymbol{v}^{\boldsymbol{T}} - \boldsymbol{F}_{rPT}^*(\boldsymbol{F}_{rPT}^*)^T \succcurlyeq \boldsymbol{0}$.

First, we show that $\boldsymbol{v}$ is parallel to $\boldsymbol{F}_{rPT}^*$ , i.e., there exists a $\rho \in \mathbb{R}$, such that $\boldsymbol{v} = \rho \boldsymbol{F}_{rPT}^*$.

Suppose that for all $\rho \in \mathbb{R}$, $\boldsymbol{v} \neq \rho \boldsymbol{F}_{rPT}^*$t, then there always exists a vector $\boldsymbol{a}$, such that $\boldsymbol{a}^{\boldsymbol{T}}\boldsymbol{v} = 0$ and $\boldsymbol{a}^{\boldsymbol{T}}\boldsymbol{F}_{rPT}^* \neq 0$. Using the same $\boldsymbol{a}$, we have:

$$\boldsymbol{a}^{\boldsymbol{T}}\left(\boldsymbol{v}\boldsymbol{v}^{\boldsymbol{T}} - \boldsymbol{F}_{rPT}^*(\boldsymbol{F}_{rPT}^*)^T\right)\boldsymbol{a} = \boldsymbol{a}^{\boldsymbol{T}}\boldsymbol{v}\boldsymbol{v}^{\boldsymbol{T}}\boldsymbol{a} - \boldsymbol{a}^{\boldsymbol{T}}\boldsymbol{F}_{rPT}^*(\boldsymbol{F}_{rPT}^*)^T\boldsymbol{a} = -\left(\boldsymbol{a}^{\boldsymbol{T}}\boldsymbol{F}_{rPT}^*\right)^{\boldsymbol{2}}$$

which contradicts with $\boldsymbol{v}\boldsymbol{v}^{\boldsymbol{T}} - \boldsymbol{F}_{rPT}^*(\boldsymbol{F}_{rPT}^*)^T \succcurlyeq 0$. Therefore, $\boldsymbol{v} = \rho \boldsymbol{F}_{rPT}^*$.

Then, we prove $\rho^2 = 1$.

(ii) Substitute $\boldsymbol{v} = \rho \boldsymbol{F}_{rPT}^*$ into constraint $(3h'')$, we have:

$$\boldsymbol{F}_{rPT}^*(l)\sum_{j=1}^{6N_1}\boldsymbol{A}^*(i_l,j)\boldsymbol{F}_g^*(j) + \rho^2\boldsymbol{F}_{rPT}^*(l)\sum_{k\in\{1,\dots,3N_{PT}\}}\boldsymbol{A}^*\big(m(l),m(k)\big)\boldsymbol{F}_{rPT}^*(k) = 0.$$

which is equivalent to

$$(1-\rho^2)\boldsymbol{F}_{rPT}^*(l)\sum_{j=1}^{6N_1}\boldsymbol{A}^*(m(l),j)\boldsymbol{F}_g^*(j) + \rho^2\boldsymbol{F}_{rPT}^*(l)\left(\sum_{j=1}^{6N_1}\boldsymbol{A}^*(m(l),j)\boldsymbol{F}_g^*(j) + \sum_{k\in\{1,\dots,3N_{PT}\}}\boldsymbol{A}^*\big(m(l),m(k)\big)\boldsymbol{F}_{rPT}^*(k)\right) = 0.$$

Since $\sum_{j=1}^{6N_1} \boldsymbol{A}^*(m(l),j)\boldsymbol{F}_g^*(j) + \sum_{k\in\{1,\dots,3N_{PT}\}} \boldsymbol{A}^*\big(m(l),m(k)\big)\boldsymbol{F}_{rPT}^*(k) = \boldsymbol{U}^*\big(m(l)\big)$, we have

$$(1-\rho^2)\boldsymbol{F}_{rPT}^*(l)\sum_{j=1}^{6N_1} \boldsymbol{A}^*(m(l),j)\boldsymbol{F}_g^*(j) + \rho^2\boldsymbol{F}_{rPT}^*(l)\boldsymbol{U}^*\big(m(l)\big) = 0.$$

Since $\boldsymbol{F}_{rPT}^*(l)\boldsymbol{U}^*\big(m(l)\big) = 0$, we have

$$(1-\rho^2)\boldsymbol{F}_{rPT}^*(l)\sum_{j=1}^{6N_1} \boldsymbol{A}^*(m(l),j)\boldsymbol{F}_g^*(j) = \rho^2\boldsymbol{F}_{rPT}^*(l)\boldsymbol{U}^*\big(m(l)\big) = 0.$$

According to the assumption that there is $l \in \{1,\dots,3N_{PT}\}$ such that $\boldsymbol{F}_{rPT}^*(l)\sum_{j=1}^{6N_1} \boldsymbol{A}^*(m(l),j)\boldsymbol{F}_g^*(j) \neq 0$, we conclude that $1-\rho^2 = 0$. Finally, the proof is complete.